\documentclass{article}

\usepackage{arxiv}

\usepackage[utf8]{inputenc} 
\usepackage[T1]{fontenc}   
\usepackage{hyperref}      
\usepackage{url}         
\usepackage{booktabs}    
\usepackage{amsfonts}  
\usepackage{nicefrac}    
\usepackage{microtype}   
\usepackage{lipsum}
\usepackage{changepage}
\usepackage{graphicx}
\usepackage{amsmath}
\usepackage{amsfonts}
\usepackage{amssymb}
\usepackage{makeidx}
\usepackage{algorithm,algorithmic}
\usepackage[switch]{lineno}
\usepackage{xcolor}
\usepackage{bm}
\usepackage[justification=centering]{caption}
\usepackage{amsfonts}
\usepackage[bbgreekl]{mathbbol}
\DeclareSymbolFontAlphabet{\mathbbm}{bbold}
\DeclareSymbolFontAlphabet{\mathbb}{AMSb}
\usepackage{amsthm}

\title{A Note on Linear Quadratic Regulator\\ and Kalman Filter}

\author{ Midhun T. Augustine\\
  Department of Electrical Engineering\\
  Indian Institute of Technology Delhi, India \\
  midhunta30@gmail.com\\
  \vspace{.05cm}\\
  Date of initial version: 10 - 02 - 2021\\ Date of current version: 30 - 08 - 2023\\ \vspace{.05cm}}

\begin{document}
\maketitle

\begin{abstract}
Two central problems in modern control theory are the controller design problem: which deals with designing a control law for the dynamical system, and the state estimation problem (observer design problem): which deals with computing an estimate of the states of the dynamical system. The Linear Quadratic Regulator (LQR) and Kalman Filter (KF) solves these problems respectively for linear dynamical systems in an optimal manner, i.e., LQR is an optimal state feedback controller and KF  is an optimal state estimator. In this note, we will be discussing the basic concepts, derivation, steady-state analysis, and numerical implementation of the LQR and KF.
\end{abstract}


\section{Introduction}
Optimal control is one of the modern control design approaches that integrates optimization to controller design, i.e., in optimal control, the control law is designed for minimizing a performance measure or cost function of the system, and the control design problem can be formulated as an optimization problem.  LQR is one of the most popular optimal control approaches, which deals with finding the optimal control input that drives the state $\textbf{x}$ of a linear system to a desired operating point $\textbf{x}_{r}$ (regulating the state at $\textbf{x}_{r}$) while minimizing a quadratic performance measure.
\par Similarly, the optimal estimation problem deals with computing an estimate of the unknown variable or vector that minimizes a performance measure. In practice, the states of a dynamical system can be considered to be unknown or partially known, and the KF deals with finding an optimal estimate of the states of the linear system, in which the performance measure is chosen as a quadratic function of the estimation error. The objective of this note is to give an intuitive introduction to the LQR and KF from a theoretical point of view. 

\par \textit{Notations:} $\mathbb{R}^{n}$ stands for $n$ - dimensional Euclidean space and $\mathbb{R}^{m \times n}$ refers to the space of $m \times n$ real matrices. $\mathbb{N},\mathbb{Z}$ and $\mathbb{C}$ denotes the set of natural numbers, integers and complex numbers respectively. 
Matrices and vectors are represented by boldface letters ($\textbf{A},\textbf{a}$),  scalars by normal font ($A,a$), and sets by blackboard bold font ($\mathbb{A},\mathbb{B},\mathbb{C}$). The notation $\textbf{P}>0$ indicates that $\textbf{P}$ is a real  positive definite matrix and $\textbf{P}\geq0$ refers to positive semidefinite matrices. Finally, $\textbf{I}$  is used to represent  identity matrix.

\section{Linear Quadratic Regulator}
 For this discussion, we are considering the discrete-time linear  time-varying (LTV) system defined by 
 \begin{equation}
\textbf{x}_{k+1}=\textbf{A}_{k}\textbf{x}_{k}+\textbf{B}_{k}\textbf{u}_{k}
 \end{equation}
where $ k\in \{0,1,..,N-1 \} \subseteq \mathbb{Z}^{+}$ is the discrete time instant, $N$ is the length of the time horizon and $\mathbb{Z}^{+}$ is the set of non-negative integers, $\textbf{x}_{k}\in \mathbb{R}^{n},$ $\textbf{u}_{k}\in \mathbb{R}^{m}$  are the state vector, input vector  and $\textbf{A}_{k}\in \mathbb{R}^{n \times n},$ $\textbf{B}_{k}\in \mathbb{R}^{n \times m}$  are the system matrix, input matrix respectively. 
 We can start the discussion with state feedback controller design using eigenvalue placement (also known as pole placement) for the discrete-time linear time-invariant (LTI) system, for which $\textbf{A}_{k}=\textbf{A},$ $\textbf{B}_{k}=\textbf{B}.$ Also, we assume
 all  states are available for feedback, i.e., $\textbf{x}_{k}$ is known.
 Consider a stabilization problem in which the objective is to drive the state to zero from any initial state, i.e., $\textbf{x}_{r}=\textbf{0}$ (in general,  regulation problem with nonzero reference: $\textbf{x}_{r}\neq \textbf{0}$
can be transferred to stabilization problem by using the coordinate transformation: $\bar{\textbf{x}}=\textbf{x}-\textbf{x}_{r},$ then stabilizing $\bar{\textbf{x}}$). 
 Now, for the linear state feedback:
\begin{equation}
\textbf{u}_{k}=-\textbf{K}\textbf{x}_{k}
\end{equation}
the state equation (1) becomes
\begin{equation}
  \textbf{x}_{k+1}=[\textbf{A}-\textbf{B}\textbf{K}]\textbf{x}_{k}
  \end{equation}
where $\textbf{A}-\textbf{B}\textbf{K}$ is the closed loop system matrix.
 We have studied in basic control theory that, if the system $[\textbf{A},\textbf{B}]$ is controllable (i.e., the rank of the matrix $\left[\begin{matrix} \textbf{B} & \textbf{A}\textbf{B} & \dots & \textbf{A}^{n-1}\textbf{B}
\end{matrix}\right]$ is equal to $n$), then the feedback gain $\textbf{K}$ can be designed for placing the  eigenvalues of $\textbf{A}-\textbf{B}\textbf{K}$ to any desired locations in the complex plane, which leads to the  eigenvalue placement problem. In eigenvalue placement, one can specify the desired eigenvalues  $\mu_{i}\in \mathbb{C},i=1,2,...,n,$  based on the transient requirement and design the feedback gain $\textbf{K}$ such that the eigenvalues of $\textbf{A}-\textbf{B}\textbf{K}$ are equal to $\mu_{i}.$ 
If all the desired eigenvalues or closed-loop eigenvalues are placed within the unit disk ($|\mu_{i}|<1$), stabilization can be achieved asymptotically, i.e., $\textbf{x}_{k}\rightarrow \textbf{0}$ as $k\rightarrow \infty.$  The Mayne-Murdoch formula gives  the feedback gain (for the system $[\textbf{A},\textbf{B}]$ represented in diagonal form and eigenvalues of $\textbf{A}$ are assumed to be distinct) in terms of the open-loop eigenvalues $\lambda_{i}$ (eigenvalues of $\textbf{A}$) and closed-loop eigenvalues $\mu_{i}$ (eigenvalues of $\textbf{A}-\textbf{B}\textbf{K}$) as: 
\begin{equation}
K_{i}= \frac{1}{B_{i}}  \dfrac{\underset{j=1,..,n}{\Pi}[\lambda_{i}-\mu_{j}] }{\underset{j=1,..,n, j\neq i}{\Pi}[\lambda_{i}-\lambda_{j}]} \hspace{1cm} i=1,...,n
\end{equation}
where  $K_{i}$ and $B_{i}$ are the $i^{th}$ element of the feedback gain matrix and input matrix.
In discrete-time systems, one can place the eigenvalues closer to the origin for faster transient response, but according to (4) the feedback gain increases as the distance between open-loop and closed-loop eigenvalues increases, and thereby the control input magnitude increases. We can say that, the eigenvalue placement only considers the transient performance, and in practice, we have to include the control effort or energy required, in the design problem. It would be desirable to have a faster transient response with lesser control effort, and one way to achieve this is to consider a performance measure or cost function as a weighted sum of states and control inputs over the time horizon, and design the control law for minimizing the cost function. This leads to the LQR in which  the cost function is chosen as  a quadratic sum of the states  and control inputs:
\begin{equation}
J=\textbf{x}_{N}^{T}\textbf{Q}_{N}\textbf{x}_{N}+ \sum_{k=0}^{N-1}\textbf{x}_{k}^{T}\textbf{Q}_{k}\textbf{x}_{k}+\textbf{u}_{k}^{T}\textbf{R}_{k}\textbf{u}_{k} 
\end{equation}
where $\textbf{Q}_{k} \in \mathbb{R}^{n \times n},$ $\textbf{R}_{k}  \in \mathbb{R}^{m \times m}$ are the weighting matrices used for relatively weighting the states and control inputs and to be chosen such that $\textbf{Q}_{k} \geq 0,$ $\textbf{R}_{k}>0$. Now the LQR problem can be defined as
\par \textbf{LQR Problem:} Find the optimal control input sequence $\textbf{u}_{k}, k=0,1,...,N-1,$ for the linear system (1) which minimizes the quadratic cost function (5).
\par The popoularity of LQR is mainly because of,\\
\begin{tabular}{@{}l@{\ }l}
	\hspace{.5cm}	i. & It gives a linear state feedback control law: $\textbf{u}_{k}=-\textbf{K}_{k}\textbf{x}_{k}.$ \\
	\hspace{.5cm}	ii. & Considers both the transient performance and control effort in designing the control law. \\
	\hspace{.5cm}	iii. & It can be easily applied to LTV systems, whereas eigenvalue placement of LTV systems is difficult. \\
\end{tabular}\newline
\subsection{LQR solution: The dynamic programming approach}
Here we are using the dynamic programming approach for solving the LQR problem. In dynamic programming, the optimal control problem is solved recursively in time using the idea of the cost-to-go function $V_{k},$ where 
 $V_{k}$ is the cost accumulated from the $k^{th}$ instant till the end, and for the quadratic cost function (5) we define $V_{k}$ as
\begin{equation}
V_{N}=\textbf{x}_{N}^{T}\textbf{Q}_{N}\textbf{x}_{N}, \hspace{1cm} V_{k}=\textbf{x}_{k}^{T}\textbf{Q}_{k}\textbf{x}_{k}+\textbf{u}_{k}^{T}\textbf{R}_{k}\textbf{u}_{k}+V_{k+1}. \hspace{0.5cm} k=N-1,...,1,0.
\end{equation}
From (5) and (6) we obtain  $J=V_{0},$
and we can solve the optimal control problem by recursively minimizing $V_{k}$ backwards, i.e., computing the optimal cost-to-go functions (also known as the value functions): $V_{N-1}^{*},...,V_{1}^{*},V_{0}^{*}$ from which the optimal cost is obtained as $J^{*}=V_{0}^{*}.$ 
For linear systems with quadratic cost, the optimal value function will be quadratic, and can be represented as
\begin{equation}
    V_{k}^{*}=\underset{\textbf{u}_{k}}{min} \hspace{.1cm} \textbf{x}_{k}^{T}\textbf{Q}_{k}\textbf{x}_{k}+\textbf{u}_{k}^{T}\textbf{R}_{k}\textbf{u}_{k}+V_{k+1}^{*}=\textbf{x}_{k}^{T}\textbf{P}_{k}\textbf{x}_{k}
\end{equation}
where $\textbf{P}_{k} \geq 0$ is known as the Riccati matrix which is to be determined, for that, we substitute (1)
in (7), results in:
\begin{equation}
V_{k}^{*}=\underset{\textbf{u}_{k}}{min} \hspace{.1cm} \textbf{x}_{k}^{T}\textbf{Q}_{k}\textbf{x}_{k}+\textbf{u}_{k}^{T}\textbf{R}_{k}\textbf{u}_{k}+[{\textbf{A}}_{k}\textbf{x}_{k}+{\textbf{B}}_{k}\textbf{u}_{k}]^{T}\textbf{P}_{k+1}[{\textbf{A}}_{k}\textbf{x}_{k}+{\textbf{B}}_{k}\textbf{u}_{k}].
\end{equation}
From which the optimal control input is obtained by equating the gradient of the cost-to-go function with respect to the control input to zero, i.e., ${\nabla V_{k}}_{\textbf{u}_{k}}=\frac{\partial V_{k}}{\partial  \textbf{u}_{k}}= \left[\begin{matrix} \frac{\partial V_{k}}{\partial  {u}_{1_k}} & \frac{\partial V_{k}}{\partial  {u}_{2_k}} & \dots & \frac{\partial V_{k}}{\partial  {u}_{m_k}}
\end{matrix}\right]^{T}=\textbf{0}$ (the first-order necessary condition for optimality which is also sufficient in this case, since $V_{k}$ is a convex function), which results in:
\begin{equation}
    \begin{aligned}
    &   2\textbf{R}_{k}\textbf{u}_{k}+2{\textbf{B}}_{k}^{T}\textbf{P}_{k+1}[{\textbf{A}}_{k}\textbf{x}_{k}+{\textbf{B}}_{k}\textbf{u}_{k}]=\textbf{0}\\
    \implies&\textbf{u}_{k}=\textbf{u}_{k}^{*}=-[\textbf{R}_{k}+{\textbf{B}}_{k}^{T}\textbf{P}_{k+1}{\textbf{B}}_{k}]^{-1}{\textbf{B}}_{k}^{T}\textbf{P}_{k+1}{\textbf{A}}_{k}\textbf{x}_{k}=-\textbf{K}_{k}\textbf{x}_{k}
\end{aligned}
\end{equation}
where 
\begin{equation}
\textbf{K}_{k}=[\textbf{R}_{k}+{\textbf{B}}_{k}^{T}\textbf{P}_{k+1}{\textbf{B}}_{k}]^{-1}{\textbf{B}}_{k}^{T}\textbf{P}_{k+1}{\textbf{A}}_{k}
\end{equation}
is the optimal feedback gain. Thus, we have the optimal control input as linear state feedback with time-varying gain.  Since,    $\textbf{R}_{k}>0$ and  $\textbf{P}_{k}\geq 0,$ we have $\textbf{R}_{k}+{\textbf{B}}_{k}^{T}\textbf{P}_{k+1}{\textbf{B}}_{k}$ will be positive definite, hence invertible.  This guarantees the existence and uniqueness of the feedback gains at each time instant.
Now, by substituting $\textbf{u}_{k}^{*}$ instead of $\textbf{u}_{k}$  in (8), we obtain the  optimal cost-to-go function or  value function as a quadratic function: 
\begin{equation}
\begin{aligned}
V_{k}^{*}=&\textbf{x}_{k}^{T}\textbf{Q}_{k}\textbf{x}_{k}+[-\textbf{K}_{k}\textbf{x}_{k}]^{T}\textbf{R}_{k}[-\textbf{K}_{k}\textbf{x}_{k}]+\big[[{\textbf{A}}_{k}-{\textbf{B}}_{k}{\textbf{K}}_{k}]\textbf{x}_{k}\big]^{T}\textbf{P}_{k+1}\big[[{\textbf{A}}_{k}-{\textbf{B}}_{k}{\textbf{K}}_{k}]\textbf{x}_{k}\big]\\
=&\textbf{x}_{k}^{T}\big[\textbf{Q}_{k}+\textbf{K}_{k}^{T}\textbf{R}_{k} \textbf{K}_{k}+ [\textbf{A}_{k}- \textbf{B}_{k}  \textbf{K}_{k}]^{T}\textbf{P}_{k+1}[\textbf{A}_{k}- \textbf{B}_{k}  \textbf{K}_{k}]\big]\textbf{x}_{k}=\textbf{x}_{k}^{T}\textbf{P}_{k}\textbf{x}_{k}
\end{aligned}
\end{equation}
where 
\begin{equation}
\textbf{P}_{k}=\textbf{Q}_{k}+\textbf{K}_{k}^{T}\textbf{R}_{k} \textbf{K}_{k}+ [\textbf{A}_{k}- \textbf{B}_{k}  \textbf{K}_{k}]^{T}\textbf{P}_{k+1}[\textbf{A}_{k}- \textbf{B}_{k}  \textbf{K}_{k}]
\end{equation}
is known as the Difference Riccati  Equation (DRE)  using which the Riccati matrix is computed recursively from a terminal Riccati matrix $\textbf{P}_{N}=\textbf{Q}_{N}$. For the LQR algorithm,
we define the sets $\mathbb{A}= \{\textbf{A}_{0},\textbf{A}_{1},...,\textbf{A}_{N-1} \},\mathbb{B}= \{\textbf{B}_{0},\textbf{B}_{1},...,\textbf{B}_{N-1}\}, \mathbb{Q}= \{\textbf{Q}_{0},\textbf{Q}_{1},...,\textbf{Q}_{N}\}, \mathbb{R}= \{\textbf{R}_{0},\textbf{R}_{1},...,\textbf{R}_{N-1}\},$ $\mathbb{K}= \{\textbf{K}_{0},\textbf{K}_{1},...,\textbf{K}_{N-1}\},$ and, let $[\mathbb{A}]_{i}$ denotes the $i^{ith}$ element of the set $\mathbb{A}.$
Now, the algorithm for computing $\mathbb{K}$ is given below: 
\begin{algorithm}[H]
 \small
	\begin{algorithmic}[1] 
	
	\STATE Require $\mathbb{A},\mathbb{B},\mathbb{Q},\mathbb{R}$
		\STATE Initialize $\textbf{P}_{N}=[\mathbb{Q}]_{N+1}$   
		\FOR  {$k= N-1~to~ 0 $}
		\STATE ${\textbf{A}}_{{k}}= [\mathbb{A}]_{k+1},$ ${\textbf{B}}_{{k}}= [\mathbb{B}]_{k+1},$ ${\textbf{Q}}_{{k}}= [\mathbb{Q}]_{k+1},$ ${\textbf{R}}_{{k}}= [\mathbb{R}]_{k+1}$
		\STATE  $\textbf{K}_{k}=[\textbf{R}_{k}+{\textbf{B}}_{k}^{T}\textbf{P}_{k+1}{\textbf{B}}_{k}]^{-1}{\textbf{B}}_{k}^{T}\textbf{P}_{k+1}{\textbf{A}}_{k}$ and $[\mathbb{K}]_{k+1}=\textbf{K}_{k}$
		\STATE $\textbf{P}_{k}=\textbf{Q}_{k}+\textbf{K}_{k}^{T}\textbf{R}_{k} \textbf{K}_{k}+ [\textbf{A}_{k}- \textbf{B}_{k}  \textbf{K}_{k}]^{T}\textbf{P}_{k+1}[\textbf{A}_{k}- \textbf{B}_{k}  \textbf{K}_{k}]$
		\ENDFOR
		\STATE $J^{*}=V_{0}^{*}=\textbf{x}_{0}^{T}\textbf{P}_{0} \textbf{x}_{0}$
		\STATE Return $\mathbb{K}$
	\end{algorithmic}
	\caption{(LQR: Backward  recursion)}
\end{algorithm}
Once the set $\mathbb{K}$ is obtained, the optimal control input  can be computed at each time instant, and the algorithm for the forward simulation of the system with optimal state feedback is given below:
\begin{algorithm}[H]
 \small
	\begin{algorithmic}[1] 
	
	\STATE Require $\mathbb{A},\mathbb{B}$ and $\mathbb{K}$
		\STATE Initialize $\textbf{x}_{0}$   
		\FOR  {$k= 0~to~ N-1 $}
		\STATE $\textbf{A}_{k}=[\mathbb{A}]_{k+1},$ $\textbf{B}_{k}=[\mathbb{B}]_{k+1},$ $\textbf{K}_{k}=[\mathbb{K}]_{k+1}$ \STATE $\textbf{u}_{k}=-\textbf{K}_{k} \textbf{x}_{k}$
		\STATE $\textbf{x}_{k+1}=\textbf{A}_{k} \textbf{x}_{k}+\textbf{B}_{k} \textbf{u}_{k}$
		\ENDFOR
	\end{algorithmic}
	\caption{(Forward simulation)}
\end{algorithm}
\subsection{LQR in steady state}
We have seen that for LTV systems, the optimal feedback gain will be time-varying, since $\textbf{A}_{k},\textbf{B}_{k},\textbf{Q}_{k},\textbf{R}_{k}$ is time-varying, and in general, we cannot say much about the convergence of $\textbf{P}_{k}$ and $\textbf{K}_{k}$. But, for LTI systems we have $\textbf{A}_{k}=\textbf{A},\textbf{B}_{k}=\textbf{B},\textbf{Q}_{k}=\textbf{Q},\textbf{R}_{k}=\textbf{R}$ and the only time-varying element in the  feedback gain (10) will be $\textbf{P}_{k+1}.$
Therefore, if $\textbf{P}_{k}$ converges to some steady state matrix $\textbf{P},$ the feedback gain $\textbf{K}_{k}$ also converges to some steady state gain matrix \textbf{K}, and we have the following result regarding the convergence:
\par For LTI systems,
if $[\textbf{A},\textbf{B}]$ is controllable and $\textbf{Q}>0,$ the DRE (12)  converges to a unique positive definite solution $\textbf{P}$ of the Algebraic Riccatti Equation (ARE): \begin{equation}
  \textbf{P}=\textbf{Q}+\textbf{K}^{T}\textbf{R} \textbf{K}+[\textbf{A}- \textbf{B}  \textbf{K}]^{T}\textbf{P}[\textbf{A}- \textbf{B} \textbf{K}] 
\end{equation}
which results in the unique feedback gain:
\begin{equation}
  \textbf{K}=[ \textbf{R}+\textbf{B}^{T}\textbf{P}\textbf{B} ]^{-1} \textbf{B}^{T}\textbf{P}\textbf{A}  
\end{equation}
such that all the eigenvalues of $\textbf{A}-\textbf{B}\textbf{K}$ lies inside the unit disk. We can also have a less restrictive condition: $[\textbf{A},\textbf{B}]$ is stabilizable, $\textbf{Q}\geq 0$ and $[\textbf{A},\textbf{Q}^{\frac{1}{2}}]$ is observable or detectable, which is sufficient for convergence.
\par One important observation is that for the states, the transient period starts from $0,$ and let $k_{\textbf{x}}$ be the time instant at which the states reaches its steady-state values (or within some $\epsilon$ - neighbourhood). Similarly, for the feedback gains the transient period starts from $N-1,$ and let $k_{\textbf{K}}$ be the time instant at which the feedback gain reaches its final values, i.e., for $0\leq k\leq k_{\textbf{K}}$ the optimal feedback gain will be a fixed matrix, and for $k_{\textbf{K}}<k\leq N-1$
the optimal feedback gain will be time-varying. Now, if we have the condition: 
\begin{equation}
    k_{\textbf{x}}<k_{\textbf{K}}
\end{equation}
then the optimal feedback gains during the transient period of the states will be a fixed matrix. We have, only the values of the feedback gains in the transient period of the states will affect the cost. Therefore, we can replace the time-varying gain $\textbf{K}_{k}$ by a fixed gain $\textbf{K}$ without affecting the cost, if condition (15) is satisfied. In this case, we need to store only one feedback gain $\textbf{K},$ which is obtained from (14) instead of storing the set $\mathbb{K},$ which contain $N$ matrices.
For a better understanding of this concept, we consider an LTI system with dynamics:
\begin{equation}
    \textbf{A}=\left[\begin{matrix} 0.5 &0 \\-1 & 1.5
\end{matrix}\right] \hspace{1cm} \textbf{B}=\left[\begin{matrix} 0.5 \\0.1
\end{matrix}\right]
\end{equation}
for which the sequence of feedback gains and Riccati matrices are computed using algorithm 1, and the plot of the states, control input, feedback gains, and diagonal elements of the Riccati matrix are given in Fig. 1, for $\textbf{Q}=\textbf{I}_{2},\textbf{R}=1$ and $\textbf{x}_{0}=\left[\begin{matrix} 10 & 5
\end{matrix}\right]^{T}.$ In this example the feedback gain converges to a fixed matrix $\textbf{K}=\left[\begin{matrix} 2.73 &-2.75
\end{matrix}\right]$ and, let $J(\textbf{K})$ denotes the cost for the system (16) with this fixed gain matrix, and $J(\textbf{K}_{k})$ denotes the cost for the optimal feedback gain. Now, for $N=5,$ the optimal feedback gain is time-varying over the transient period of the states (see Fig. 1(a)). Therefore,  replacing the time-varying gain with a fixed gain results in  loss of optimality, and for this example we obtained $J(\textbf{K}_{k})=422.13,$  ${J}(\textbf{K})=432.17$ in which $J(\textbf{K}_{k})<J(\textbf{K})$. But, for $N=50,$ the time horizon is sufficiently large that results in the feedback gain to converge, and the optimal feedback gain is constant over the transient period of the states, i.e., condition (15) is satisfied (see Fig. 1(b)). Therefore,  we can replace the time-varying gain by $\textbf{K}$ without affecting the cost, and for this example we obtained $J(\textbf{K}_{k})=433.25,$ $J(\textbf{K})=433.25,$ in which $J(\textbf{K}_{k})=J(\textbf{K})$.

Consequently, for infinite horizon problems, if the states and feedback gains reaches to their steady state values in a finite number of time instants, then we have clearly  $k_{\textbf{x}}<k_{\textbf{K}},$ since $k_{\textbf{x}}$ ill be a finite value and $k_{\textbf{K}}=\infty-finite \hspace{0.1cm} value,$ will be infinity. And for controllable or stabilizable linear systems, the states and optimal feedback gains converges to their steady state values (with an $\epsilon$ - error) in finite number of time instants.
Therefore, for the infinite horizon LQR of  controllable or stabilizable linear systems, the optimal control law will be a time-invariant feedback: $\textbf{u}_{k}^{*}=-\textbf{K}\textbf{x}_{k}$ where $\textbf{K}$ is  obtained from  (14).

\begin{figure} [H]	
	\begin{center}
		\includegraphics [scale=.63] {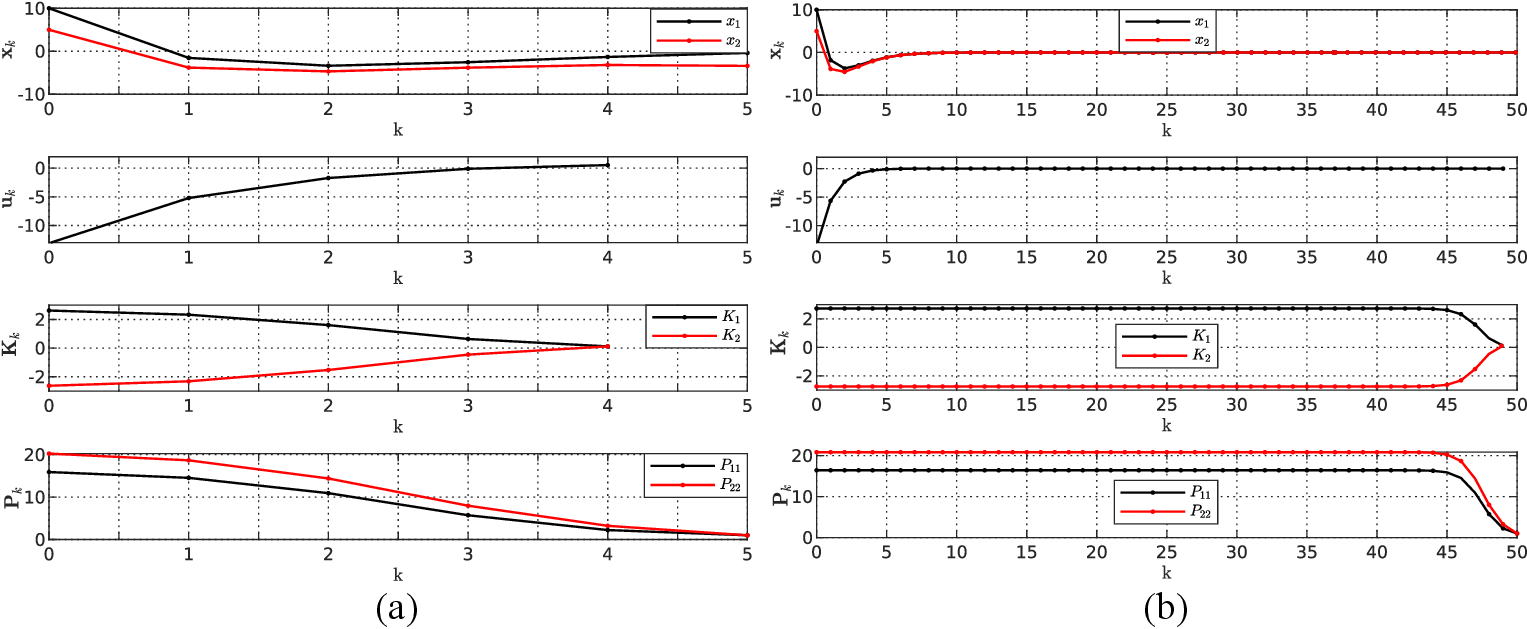}
		\caption{{\footnotesize Response of LTI system with LQR:  (a) $N=5,$ (b) $N=50$.}}	
	\end{center}
\end{figure}

\section{Kalman Estimator}
In the previous section on LQR design, we have assumed that all the states are available for feedback. But measuring all the states using sensors is impractical, and usually, we measure a reduced number of variables, say $p\leq n$ using sensors, which are called the output variables. Also, in practice, the system may be affected by uncertain inputs such as disturbances and measurement noise, therefore, we have a more general model of the discrete-time linear system:
\begin{equation}
 \begin{aligned}
 &\textbf{x}_{k+1}=\textbf{A}_{k}\textbf{x}_{k}+\textbf{B}_{k}\textbf{u}_{k}+\textbf{d}_{k}\\
 &\textbf{y}_{k}=\textbf{C}_{k}\textbf{x}_{k}+\textbf{v}_{k}
  \end{aligned}
 \end{equation}
where $\textbf{x}_{k}\in \mathbb{R}^{n},$ $\textbf{u}_{k}\in \mathbb{R}^{m},$ $\textbf{y}_{k}\in \mathbb{R}^{p},$ $\textbf{d}_{k}\in \mathbb{R}^{n},$ $\textbf{v}_{k}\in \mathbb{R}^{p}$ are the state vector, input vector, output vector, disturbance vector, noise vector  and $\textbf{A}_{k}\in \mathbb{R}^{n \times n},$ $\textbf{B}_{k}\in \mathbb{R}^{n \times m},$ $\textbf{C}_{k}\in \mathbb{R}^{p \times n}$ are the system matrix, input matrix, output matrix respectively. The system (17) is also known as the discrete-time stochastic linear system,  in which the disturbance $\textbf{d}_{k},$ noise $\textbf{v}_{k}$ and initial state $\textbf{x}_{0}$ are random vectors which makes the state $\textbf{x}_{k}$ and output $\textbf{y}_{k}$ stochastic.
Then, the estimation problem is to compute an estimate of the  state $\textbf{x}_{k}$ denoted by $\hat{\textbf{x}}_{k}$
using the available information, which includes the system model $[\textbf{A}_{k},\textbf{B}_{k},\textbf{C}_{k}],$ input $\textbf{u}_{k},$ output $\textbf{y}_{k}$ and some statistical information about the random vectors $\textbf{x}_{0},\textbf{d}_{k},\textbf{v}_{k}.$
\par Let $l$ denotes the time instant upto which the output or measurement data is available, and suppose the estimate of the state: $\hat{\textbf{x}}_{k}$ at time instant $k$ is made using the measurement data upto  time instant $l,$ i.e., $\{\textbf{y}_{0},\textbf{y}_{1},...,\textbf{y}_{l} \},$ then the estimation problem can be classified as follows:\\
\begin{tabular}{@{}l@{\ }l}
	\hspace{.5cm}	i. & Prediction problem: in which $l<k$ and the estimator is called predictor.\\
	\hspace{.5cm}	ii. &Filtering problem: in which $l = k$ and the estimator is called filter.\\
		\hspace{.5cm}	iii. &Smoothing  problem: in which $l > k$ and the estimator is called smoother.
\end{tabular}\newline
In this section, we will be discussing one of the most  popular linear estimator called the Kalman estimator, which gives an optimal estimate of the state $\textbf{x}_{k}$ of the stochastic linear system (17), and for the Kalman estimator we can have the above three versions as: Kalman predictor, Kalman filter, and Kalman smoother, which are used for a wide range of applications in control systems, guidance and navigation, signal processing, economics, medical science etc. Before discussing the Kalman estimator, which is a stochastic estimator, we discuss the Luenberger observer, which is a popular deterministic estimator. And
  for the deterministic LTI system:
  \begin{equation}
 \begin{aligned}
 &\textbf{x}_{k+1}=\textbf{A}\textbf{x}_{k}+\textbf{B}\textbf{u}_{k}\\
 &\textbf{y}_{k}=\textbf{C}\textbf{x}_{k}
  \end{aligned}
 \end{equation}
  the Luenberger observer can be defined as:
 \begin{equation}
 \hat{\textbf{x}}_{k+1}=\textbf{A}\hat{\textbf{x}}_{k}+\textbf{B}\textbf{u}_{k}+\textbf{L}[\textbf{y}_{k}-\hat{\textbf{y}}_{k}].
 \end{equation}
We define the estimation error vector $\textbf{x}_{e_{k}}=\textbf{x}_{k}-\hat{\textbf{x}}_{k},$ and the error dynamics for the Luenberger observer is 
\begin{equation}
\begin{aligned}
{\textbf{x}}_{e_{k+1}}&={\textbf{x}}_{k+1}-\hat{\textbf{x}}_{k+1}=\textbf{A}\textbf{x}_{k}+\textbf{B}\textbf{u}_{k}-\textbf{A}\hat{\textbf{x}}_{k}-\textbf{B}\textbf{u}_{k}-\textbf{L}[\textbf{C}\textbf{x}_{k}-\textbf{C}\hat{\textbf{x}}_{k}]\\
&=[\textbf{A}-\textbf{L}\textbf{C}]{\textbf{x}}_{e_{k}}
 \end{aligned}
 \end{equation}
i.e., the observer error dynamics is defined by the closed-loop observer matrix $\textbf{A}-\textbf{L}\textbf{C}.$ This leads to the observer eigenvalue placement problem, in which one has to design the observer gain matrix $\textbf{L}$ to place the eigenvalues of closed-loop matrix $\textbf{A}-\textbf{L}\textbf{C}$ at desired locations. Similarly, if the system $[\textbf{A},\textbf{C}]$ is observable (i.e.,  rank of the  matrix $\left[\begin{matrix} \textbf{C}^{T} & (\textbf{C}\textbf{A})^{T} & \dots & (\textbf{C}\textbf{A}^{n-1})^{T}
\end{matrix}\right]^{T}$ is equal to $n$), then we can place all the eigenvalues of $\textbf{A}-\textbf{L}\textbf{C}$ at desired locations in the complex plane. Moreover, if all the  eigenvalues are placed within the unit disk, then ${\textbf{x}}_{e_{k}}$ converges to zero asymptotically, i.e. $\hat{\textbf{x}}_{k}$ converges to ${\textbf{x}}_{k}$.
And for faster error convergence, the closed-loop eigenvalues should be closer to the origin, which leads to larger observer gains. Now, in the presence of disturbances and noise, i.e., $\textbf{d}_{k}\neq \textbf{0},\textbf{v}_{k}\neq\textbf{0},$ we obtain the error dynamics as 
\begin{equation}
\begin{aligned}
{\textbf{x}}_{e_{k+1}}&={\textbf{x}}_{k+1}-\hat{\textbf{x}}_{k+1}=\textbf{A}\textbf{x}_{k}+\textbf{B}\textbf{u}_{k}+\textbf{d}_{k}-\textbf{A}\hat{\textbf{x}}_{k}-\textbf{B}\textbf{u}_{k}-\textbf{L}[\textbf{C}\textbf{x}_{k}+\textbf{v}_{k}-\textbf{C}\hat{\textbf{x}}_{k}]\\&=[\textbf{A}-\textbf{L}\textbf{C}]{\textbf{x}}_{e_{k}}+\textbf{d}_{k}-\textbf{L}\textbf{v}_{k}
 \end{aligned}
 \end{equation}
which indicates, for larger observer gain $\textbf{L}$ the effect of noise will be more. Therefore, for faster error convergence  the observer gain is to be large, and for noise rejection, the observer gain is to be small, i.e., we have contrasting requirements for which a Luenberger observer with fixed gain $\textbf{L}$ may not be sufficient, and  one can go for the Kalman estimator, which is an optimal state estimator with time-varying gain $\textbf{L}_{k}.$
\par For the Kalman estimator derivation, we will be using the discrete-time linear system (17) in which the disturbance $\textbf{d}_{k},$ noise $\textbf{v}_{k}$ and initial state $\textbf{x}_{0}$ are random vectors, therefore we need some concepts from statistics and probability theory for their characterization which are briefly reviewed below: 
\newpage
\textbf{Probability space:} We denote the probability space as $(\mathbb{S}, \mathbb{E}, Pr),$ which contains:\\  
\begin{tabular}{@{}l@{\ }l}
	\hspace{.5cm}	i. & Sample space $\mathbb{S}$:  which is the set of all possible outcomes.\\
	\hspace{.5cm}	ii. &  Event space $\mathbb{E}$: which contains the set of all events (i.e., the set of all subsets of $\mathbb{S}$).\\
	\hspace{.5cm}	iii. &  Probability function $Pr$: which assigns each event in $\mathbb{E}$ a real number between 0 and 1, known as the probabi-\\&lity  measure of the event, i.e., $Pr:\mathbb{E} \rightarrow [0,1].$
\end{tabular}\newline
\textbf{Random variable:} A random variable $x$ is a real-valued function defined on the  sample space $\mathbb{S},$ i.e., $x:\mathbb{S} \rightarrow \mathbb{R},$  and a random vector $\textbf{x}=\left[\begin{matrix} x_{1} & x_{2} & \dots & x_{n}
\end{matrix}\right]^{T}$ contains random variables as its elements, i.e., $\textbf{x}:\mathbb{S} \rightarrow \mathbb{R}^{n}.$
A random variable $x$ can be associated with a
probability density function (PDF): $f(x)$, which is the function that assigns each random variable, a number between 0 and 1, i.e. $f: \mathbb{R} \rightarrow [0,1]$ similarly, for a random vector $\textbf{x}$, we have $f: \mathbb{R}^{n} \rightarrow [0,1].$ For continuous random variables, the PDF is used to specify the probability of the random variable to take a value within an interval, i.e.,  $Pr(a \leq x \leq b)=\int_{a}^{b}f(x)dx.$ Similarly, for discrete random variables, the PDF, which is known by the name probability mass function (PMF), is used to specify the probability of the random variable to take a particular value, i.e., $Pr(x=a)=f(a).$\\
\textbf{Expectation:} For a random variable with probability density function $f(x)$ the expectation/mean value is a parameter that indicates the average value, and  is defined as $E(x)= \int x f(x)dx =  \sum x_{i} f(x_{i}) $ where the integration is used for continuous random variables and summation is for discrete random variables. Similarly, for a random vector $\textbf{x}=\left[\begin{matrix} x_{1} & x_{2} & \dots & x_{n}
\end{matrix}\right]^{T}$ the expectation is defined as $\textbf{E}(\textbf{x})=\left[\begin{matrix} E(x_{1}) & E(x_{2}) & \dots & E(x_{n})
\end{matrix}\right]^{T}.$\\
\textbf{Variance:} for a  random variable $x,$ we define variance: $V(x)=E([x-E(x)][x-E(x)])=E([x-E(x)]^{2})$ which indicates the measure of spread or deviation from the mean.  Similarly, for a random vector $\textbf{x},$ the variance is defined as $\textbf{V}(\textbf{x})=\textbf{E}\big([\textbf{x}-\textbf{E}(\textbf{x})][\textbf{x}-\textbf{E}(\textbf{x})]^{T}\big),$  which is known by the names: variance matrix,covariance matrix, variance-covariance matrix, dispersion matrix, etc. The dependence between two random variables can be characterised by a parameter called covariance, and for the random variables $x$ and $y$ we define the covariance: $V(x,y)=E([x-E(x)][y-E(y)]).$ Similarly,  for the random vectors    $\textbf{x}\in \mathbb{R}^{n}$ and $\textbf{y}\in \mathbb{R}^{p}$ we define the covariance matrix:
$\textbf{V}(\textbf{x},\textbf{y})=\textbf{E}\big([\textbf{x}-\textbf{E}(\textbf{x})][\textbf{y}-\textbf{E}(\textbf{y})]^{T}\big),$ and we have $\textbf{V}(\textbf{x},\textbf{y})\in \mathbb{R}^{n \times p},$ $\textbf{V}(\textbf{x},\textbf{x})=\textbf{V}(\textbf{x})\in \mathbb{R}^{n \times n}.$ Also,  $\textbf{V}(\textbf{x},\textbf{y})=\textbf{0},$ if the random vectors $\textbf{x}$ and $\textbf{y}$ are independent.\\
\textbf{Gaussian distribution (Normal distribution):} Is the most widely used probability distribution,  which is characterised  by the Gaussian probability density function:
\begin{equation}
    f(x)=\frac{1}{\sqrt{2\pi V(x)}} e^{-\frac{1}{2}\frac{[x-E(x)]^{2}}{V(x)}}
\end{equation}
and the Gaussian distribution is shown graphically in Fig. 2(a), in which,  for $x=E(x)$ we obtain $f(x)=\frac{1}{\sqrt{2\pi V(x)}},$ and for $x=E(x)\pm \sqrt{V(x)}$ we obtain $f(x)=\frac{1}{\sqrt{2\pi V(x)}} e^{-\frac{1}{2}}=\frac{1}{\sqrt{2\pi e V(x)}}.$ From which, we can say that as $V(x)$ decreases, the distribution becomes narrower (i.e., peak value increases and width decreases), and randomness in $x$ decreases (see Fig. 2(b)). Similarly, we can define the multivariate Gaussian distribution (also known as the joint Gaussian distribution),  for the random vector $\textbf{x}$ as:
\begin{equation}
    f(\textbf{x})=\frac{1}{\sqrt{(2\pi)^{n} det(\textbf{V}(\textbf{x}))}} e^{-\frac{1}{2}[\textbf{x}-\textbf{E}(\textbf{x})]^{T} \textbf{V}(\textbf{x})^{-1}[\textbf{x}-\textbf{E}(\textbf{x})]}
\end{equation}
which is also denoted as $f(x_{1},...,x_{n}),$ and for the two random variables $x,y$  the bivariate Gaussian distribution is plotted in Fig. 3.
\begin{figure} [H]
	
	\begin{center}
		\includegraphics [scale=.6] {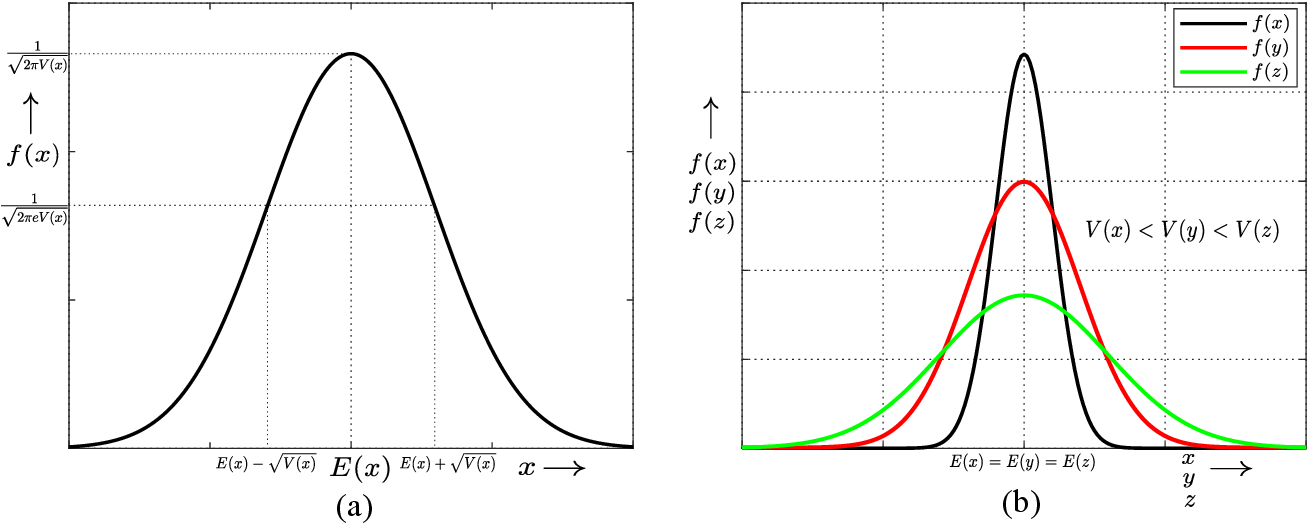}
		\caption{{\footnotesize Gaussian distribution: graphical representation}}	
	\end{center}
\end{figure} 
\textbf{Conditional probability}: 
Let $x$ and $y$ be two random variables, the conditional PDF $f(x|y)$ denotes the probability of $x$ given $y$, i.e., $f(x|y=a)$ gives the probability of $x$ given that $y$ takes the value $a,$ where $a\in \mathbb{R}$ is any scalar. For jointly Gaussian random variables $x$ and $y,$ we have $f(x|y=a)$ is Gaussian (see Fig. 3(b)). 
For, the random variables $x$ and $y$ we define the conditional expectation: $E(x|y)=\sum x_{i}f(x_{i}|y).$ Similarly, for the random vectors $\textbf{x}$ and $\textbf{y},$ we have  $\textbf{E}(\textbf{x}|\textbf{y})$ denotes the conditional expectation, and we can derive an equation for the conditional expectation, for that, we use  the following transformation: $\textbf{z}=\textbf{x}-\textbf{Ly}$
where $\textbf{L}=\textbf{V}(\textbf{x},\textbf{y})\textbf{V}(\textbf{y})^{-1}.$ Then we have: 
\begin{equation}
\begin{aligned}
   \textbf{V}(\textbf{z},\textbf{y}) =\textbf{V}(\textbf{x}-\textbf{Ly},\textbf{y})=\textbf{V}(\textbf{x},\textbf{y})
 -\textbf{L}\textbf{V}(\textbf{y})=\textbf{V}(\textbf{x},\textbf{y})-\textbf{V}(\textbf{x},\textbf{y})\textbf{V}(\textbf{y})^{-1}\textbf{V}(\textbf{y})=\textbf{0}
 \end{aligned}
 \end{equation}
which implies $\textbf{z}$ and $\textbf{y}$ are independent.  Using this we obtain: 
\begin{equation}
\begin{aligned}
   \textbf{E}(\textbf{x}|\textbf{y}) &= \textbf{E}(\textbf{z}+\textbf{L}\textbf{y}|\textbf{y})=\textbf{E}(\textbf{z}|\textbf{y})+\textbf{E}(\textbf{L}\textbf{y}|\textbf{y})\\
   &=\textbf{E}(\textbf{z})+\textbf{L}\textbf{y}=\textbf{E}(\textbf{x})-\textbf{L}\textbf{E}(\textbf{y})+  \textbf{L}\textbf{y} \hspace{1cm} \because \hspace{.1cm} \textbf{z}  \hspace{.1cm}and  \hspace{.1cm} \textbf{y}  \hspace{.1cm} are \hspace{.1cm} independent \hspace{.1cm} \textbf{E}(\textbf{z}|\textbf{y})=\textbf{E}(\textbf{z}).\\
   &=\textbf{E}(\textbf{x})+\textbf{L}[\textbf{y}-\textbf{E}(\textbf{y})]=\textbf{E}(\textbf{x})+\textbf{V}(\textbf{x},\textbf{y})\textbf{V}(\textbf{y})^{-1}[\textbf{y}-\textbf{E}(\textbf{y})].
 \end{aligned}
 \end{equation}
Similarly, we can obtain an expression for the conditional variance: 
 \begin{equation}
  \textbf{V}(\textbf{x}|\textbf{y})=\textbf{V}(\textbf{x})-\textbf{L} \textbf{V}(\textbf{y},\textbf{x})= \textbf{V}(\textbf{x})-\textbf{V}(\textbf{x},\textbf{y}) \textbf{V}(\textbf{y})^{-1} \textbf{V}(\textbf{y},\textbf{x}).
 \end{equation}
 Note that, here the selection of $\textbf{L}=\textbf{V}(\textbf{x},\textbf{y})\textbf{V}(\textbf{y})^{-1}$ is the significant part which results in the conditional expectation $\textbf{E}(\textbf{x}|\textbf{y})$ with minimum mean square error.
\begin{figure} [H]
	
	\begin{center}
		\includegraphics [scale=.7] {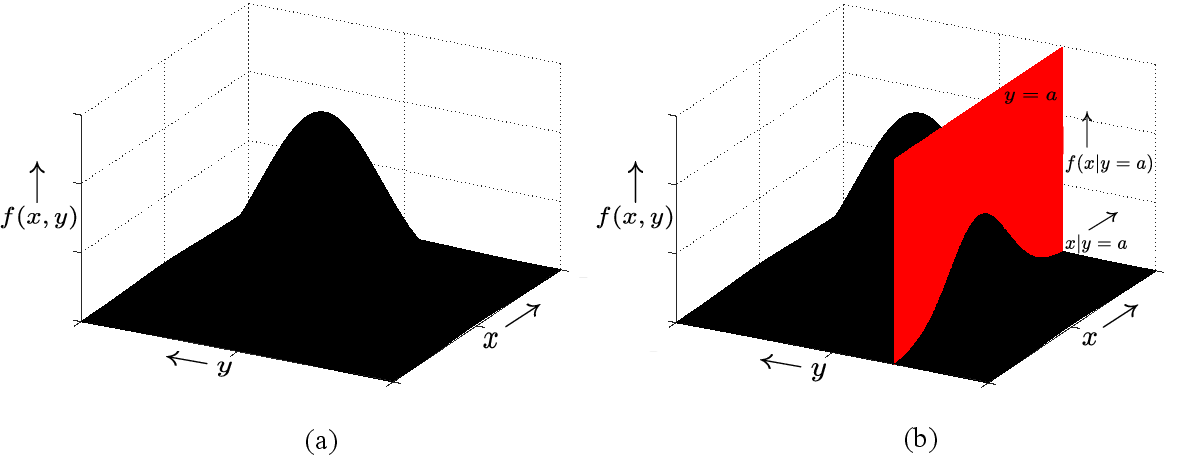}
		\caption{{\footnotesize Multivariate Gaussian distribution: graphical representation}}
	\end{center}
\end{figure} 
For the  Kalman estimator derivation and analysis, we assume that the disturbance, noise and initial state vector as Gaussian with mean $\textbf{E}({\textbf{d}}_{k})=\textbf{0},\textbf{E}( {\textbf{v}}_{k})=\textbf{0}$ and $\textbf{E}( {\textbf{x}}_{0})$ can be nonzero which is assumed to be known. Also, the vectors ${\textbf{x}}_{k},{\textbf{d}}_{k},{\textbf{v}}_{k}$ are assumed to be independent or uncorrelated i.e., $\textbf{V}(\textbf{x}_{k},\textbf{d}_{k})=\textbf{0},\textbf{V}(\textbf{x}_{k},\textbf{v}_{k})=\textbf{0},\textbf{V}(\textbf{d}_{k},\textbf{v}_{k})=\textbf{0}.$ Because of the randomness associated with ${\textbf{x}}_{0},{\textbf{d}}_{k},{\textbf{v}}_{k},$ the state vector $\textbf{x}_{k}$ and output vector $\textbf{y}_{k}$  will be a random vectors, i.e., $\textbf{x}_{k}:\mathbb{Z}^{+}\rightarrow \mathbb{R}^{n}$ and $\textbf{y}_{k}:\mathbb{Z}^{+}\rightarrow \mathbb{R}^{p}$.
The Kalman estimator gives an estimate of the states in terms of it's expectation and variance, given the  measurement data and the system model. Up next, we discuss the following versions of the Kalman estimator: predictor, filter, and smoother and there are a number of approaches for deriving the Kalman estimator which can be broadly classified into state space (time-domain) and transfer function (frequency domain) approaches. Here we focus on the state space based derivation in which we discuss the following two approaches:\\
\begin{tabular}{@{}l@{\ }l}
	\hspace{.5cm}	i. & Optimization based approach: In this the estimator gain is computed from minimizing the cost function which \\& is chosen as a quadratic function of the estimation error.\\
	\hspace{.5cm}	ii. &  Statistical approach: in which the estimator gain and variance is computed from the conditional expectation \\&and variance equations. 
\end{tabular}
\subsection{Kalman predictor}
In Kalman predictor the measurements upto $k-1^{th}$ time instant are used for computing the estimate $\hat{\textbf{x}}_{k},$ i.e., $l=k-1.$ The predictor form is similar to the Luenberger observer and defined as
 \begin{equation}
 \hat{\textbf{x}}_{k+1}=\textbf{A}_{k}\hat{\textbf{x}}_{k}+\textbf{B}_{k}\textbf{u}_{k}+\textbf{L}_{k}[\textbf{y}_{k}-\hat{\textbf{y}}_{k}].
 \end{equation}
 In which, $\textbf{A}_{k}\hat{\textbf{x}}_{k}+\textbf{B}_{k}\textbf{u}_{k}$ is known as the prediction or forecast term, and $\textbf{y}_{k}-\hat{\textbf{y}}_{k}$ is known as the correction or innovation term.
 Note that, here we are using a varying gain $\textbf{L}_{k}$ instead of the fixed gain $\textbf{L}$. We start with the optimization based derivation in which we
 define the estimation error vector  ${\textbf{x}}_{e_{k}}={\textbf{x}}_{k}-\hat{\textbf{x}}_{k},$ and  from (17) and (27) we obtain the error dynamics as
\begin{equation}
\begin{aligned}
{\textbf{x}}_{e_{k+1}}&={\textbf{x}}_{k+1}-\hat{\textbf{x}}_{k+1}=\textbf{A}_{k}\textbf{x}_{k}+\textbf{B}_{k}\textbf{u}_{k}+\textbf{d}_{k}-\textbf{A}_{k}\hat{\textbf{x}}_{k}-\textbf{B}_{k}\textbf{u}_{k}-\textbf{L}_{k}[\textbf{C}_{k}\textbf{x}_{k}+\textbf{v}_{k}-\textbf{C}_{k}\hat{\textbf{x}}_{k}]\\&=[\textbf{A}_{k}-\textbf{L}_{k}\textbf{C}_{k}]{\textbf{x}}_{e_{k}}+\textbf{d}_{k}-\textbf{L}_{k}\textbf{v}_{k}.
 \end{aligned}
 \end{equation}
 From which we obtain the variance  of the estimation error:
 \begin{equation*}
\begin{aligned}
    \textbf{V}({\textbf{x}}_{e_{k+1}})&=\textbf{E}({\textbf{x}}_{e_{k+1}}{\textbf{x}}_{e_{k+1}}^{T})=\textbf{E}\big(\big[[\textbf{A}_{k}-\textbf{L}_{k}\textbf{C}_{k}]\textbf{x}_{e_{k}}+\textbf{d}_{k}-\textbf{L}_{k}\textbf{v}_{k}\big]\big[[\textbf{A}_{k}-\textbf{L}_{k}\textbf{C}_{k}]\textbf{x}_{e_{k}}+\textbf{d}_{k}-\textbf{L}_{k}\textbf{v}_{k}\big]^{T} \big)\\
    &=\textbf{E}\big([\textbf{A}_{k}-\textbf{L}_{k}\textbf{C}_{k}]\textbf{x}_{e_{k}}\textbf{x}_{e_{k}}^{T}[\textbf{A}_{k}-\textbf{L}_{k}\textbf{C}_{k}]^{T}+\textbf{d}_{k}\textbf{d}_{k}^{T}+\textbf{L}_{k}\textbf{v}_{k}\textbf{v}_{k}^{T}\textbf{L}_{k}^{T}\big), \hspace{.1cm} \because \textbf{E}(\textbf{x}_{e_{k}}\textbf{d}_{k}^{T}) =\textbf{E}(\textbf{x}_{e_{k}}\textbf{v}_{k}^{T}) =\textbf{E}(\textbf{d}_{k} \textbf{v}_{k}^{T})=\textbf{0} \\
    &=[\textbf{A}_{k}-\textbf{L}_{k}\textbf{C}_{k}] \textbf{V}(\textbf{x}_{e_{k}})[\textbf{A}_{k}-\textbf{L}_{k}\textbf{C}_{k}]^{T}+\textbf{V}(\textbf{d}_{k})+\textbf{L}_{k}\textbf{V}(\textbf{v}_{k})\textbf{L}_{k}^{T}.
   \\
    \end{aligned}
\end{equation*}
By defining  $\textbf{P}_{k}=\textbf{V}({\textbf{x}_{e_{k}}}),$ $\textbf{Q}_{k}=\textbf{V}({\textbf{d}_{{k}}})$ and $\textbf{R}_{k}=\textbf{V}({\textbf{v}_{{k}}})$ we obtain 
\begin{equation}
    \textbf{P}_{k+1}=[\textbf{A}_{k}-\textbf{L}_{k}\textbf{C}_{k}]\textbf{P}_{k}[\textbf{A}_{k}-\textbf{L}_{k}\textbf{C}_{k}]^{T}+\textbf{Q}_{k}+\textbf{L}_{k}\textbf{R}_{k}\textbf{L}_{k}^{T}
\end{equation}
which is  the DRE for the Kalman predictor, and  $\textbf{P}_{k}\in \mathbb{R}^{n \times n},$ $\textbf{Q}_{k}\in \mathbb{R}^{n \times n}$ and $\textbf{R}_{k}\in \mathbb{R}^{p \times p}.$ If we assume $\textbf{E}(\textbf{x}_{k})=\hat{\textbf{x}}_{k},$ then we have 
$\textbf{V}(\textbf{x}_{k})=\textbf{E}\big([\textbf{x}_{k}-\textbf{E}(\textbf{x}_{k})][\textbf{x}_{k}-\textbf{E}(\textbf{x}_{k})]^{T}\big)=\textbf{E}\big([\textbf{x}_{k}-\hat{\textbf{x}}_{k}][\textbf{x}_{k}-\hat{\textbf{x}}_{k}]^{T}\big)=\textbf{E}(\textbf{x}_{e_{k}}\textbf{x}_{e_{k}}^{T})=\textbf{V}(\textbf{x}_{e_{k}})=\textbf{P}_{k},$
i.e.,
we can characterize the state vector which is a random vector by its expectation $\textbf{E}(\textbf{x}_{k})$ and variance $\textbf{V}(\textbf{x}_{k}),$ and the Kalman estimator gives an estimate of $\textbf{E}(\textbf{x}_{k}),\textbf{V}(\textbf{x}_{k})$ as $\hat{\textbf{x}}_{k},$ $\textbf{P}_{k},$ which can be computed recursively using (27) and (29). Now, the task is to optimize the estimate, and in  Kalman predictor, the performance measure is chosen as a quadratic function of the estimation error (also known as the mean square error cost function):
\begin{equation}
\begin{aligned}
    J=E(\textbf{x}_{e_{k+1}}^{T}\textbf{x}_{e_{k+1}})=E(Trace(\textbf{x}_{e_{k+1}}\textbf{x}_{e_{k+1}}^{T}))=Trace(\textbf{E}(\textbf{x}_{e_{k+1}}\textbf{x}_{e_{k+1}}^{T}))=Trace(\textbf{P}_{k+1}).
    \end{aligned}
\end{equation}
We choose the kalman gain $\textbf{L}_{k}$ to minimize the cost $J,$ and this leads to  $\frac{\partial Trace(\textbf{P}_{k+1})}{\partial  \textbf{L}_{k}}=\textbf{0},$ which results in:
\begin{equation}
    \begin{aligned}
       & -2[\textbf{A}_{k}-\textbf{L}_{k}\textbf{C}_{k}]\textbf{P}_{k}\textbf{C}_{k}^{T}+2\textbf{L}_{k}\textbf{R}_{k}=\textbf{0}\\
        &\implies \textbf{L}_{k}=\textbf{A}_{k}\textbf{P}_{k}\textbf{C}_{k}^{T}[\textbf{C}_{k}\textbf{P}_{k}\textbf{C}_{k}^{T}+\textbf{R}_{k}]^{-1}.
    \end{aligned}
\end{equation}
We can also obtain the Kalman predictor equations using the statistical approach, in which we use conditional expectation equation (25) for obtaining $\hat{\textbf{x}}_{k}=\textbf{E}(\textbf{x}_{k}|\textbf{y}_{k-1}),$ and we have
\begin{equation}
 \begin{aligned}
 \hat{\textbf{x}}_{k+1}= \textbf{E}(\textbf{x}_{k+1}|\textbf{y}_{k})= \textbf{E}(\textbf{x}_{k+1})+\textbf{L}_{k}[\textbf{y}_{k}-\textbf{E}(\textbf{y}_{k})]=\textbf{A}_{k}\hat{\textbf{x}}_{k}+\textbf{B}_{k}\textbf{u}_{k}+\textbf{L}_{k}[\textbf{y}_{k}-\hat{\textbf{y}}_{k}].
 \end{aligned}
 \end{equation}
in which $\textbf{L}_{k}$ can be obtained as
\begin{equation}
    \begin{aligned}
        \textbf{L}_{k}&=\textbf{V}(\textbf{x}_{k+1},\textbf{y}_{k})\textbf{V}(\textbf{y}_{k})^{-1}=\textbf{V}(\textbf{A}_{k}\textbf{x}_{k}+\textbf{B}_{k}\textbf{u}_{k}+\textbf{d}_{k},\textbf{C}_{k}\textbf{x}_{k}+\textbf{v}_{k})\textbf{V}(\textbf{C}_{k}\textbf{x}_{k}+\textbf{v}_{k})^{-1}\\
        &=\textbf{E}(\textbf{A}_{k}\textbf{x}_{e_{k}}\textbf{x}_{e_{k}}^{T}\textbf{C}_{k}^{T})
       \textbf{E}(\textbf{C}_{k}\textbf{x}_{e_{k}}\textbf{x}_{e_{k}}^{T}\textbf{C}_{k}^{T}+\textbf{R}_{k})^{-1}=\textbf{A}_{k}\textbf{P}_{k}\textbf{C}_{k}^{T}[\textbf{C}_{k}\textbf{P}_{k}\textbf{C}_{k}^{T}+\textbf{R}_{k}]^{-1}.
    \end{aligned}
\end{equation}
 Similarly, we have $\textbf{P}_{k+1}=\textbf{V}(\textbf{x}_{k+1}|\textbf{y}_{k})$ which can be derived from the conditional variance equation, for that we rewrite (29) as
\begin{equation}
\begin{aligned}
  \textbf{P}_{k+1}&=\textbf{A}_{k}\textbf{P}_{k}\textbf{A}_{k}^{T}+\textbf{Q}_{k}- 2\textbf{L}_{k} \textbf{C}_{k} \textbf{P}_{k} \textbf{A}_{k}^{T} + \textbf{L}_{k}[\textbf{C}_{k} \textbf{P}_{k}\textbf{C}_{k}^{T}+ \textbf{R}_{k}] \textbf{L}_{k}^{T} \hspace{0.5cm} 
sub. \hspace{.1cm}  \textbf{L}_{k} \hspace{.1cm} from (33)\\
&=\textbf{A}_{k}\textbf{P}_{k}\textbf{A}_{k}^{T}+\textbf{Q}_{k}  -\textbf{A}_{k}\textbf{P}_{k }\textbf{C}_{k}^{T}[\textbf{C}_{k} \textbf{P}_{k}\textbf{C}_{k}^{T}+ \textbf{R}_{k}]^{-1}
\textbf{C}_{k}\textbf{P}_{k}\textbf{A}_{k}^{T}
\end{aligned}
 \end{equation}
and from the conditional variance equation (26), we obtain:
\begin{equation}
\begin{aligned}
  \textbf{P}_{k+1}&=\textbf{V}(\textbf{x}_{k+1}|\textbf{y}_{k})=  \textbf{V}(\textbf{x}_{k+1})-\textbf{L}_{k}\textbf{V}(\textbf{y}_{k},\textbf{x}_{k+1})=\textbf{E}(\textbf{A}_{k}\textbf{x}_{e_{k}}\textbf{x}_{e_{k}}^{T}\textbf{A}_{k}^{T}+\textbf{Q}_{k})-\textbf{L}_{k} \textbf{E}(\textbf{C}_{k}\textbf{x}_{e_{k}}\textbf{x}_{e_{k}}^{T}\textbf{A}_{k}^{T})\\
  &=\textbf{A}_{k}\textbf{P}_{k}\textbf{A}_{k}^{T}+\textbf{Q}_{k}  -\textbf{A}_{k}\textbf{P}_{k }\textbf{C}_{k}^{T}[\textbf{C}_{k} \textbf{P}_{k}\textbf{C}_{k}^{T}+ \textbf{R}_{k}]^{-1}
\textbf{C}_{k}\textbf{P}_{k}\textbf{A}_{k}^{T}.
\end{aligned}
 \end{equation}
For the Kalman predictor algorithm, we define the sets  $\mathbb{A}= \{\textbf{A}_{0},\textbf{A}_{1},...,\textbf{A}_{N-1} \},\mathbb{B}= \{\textbf{B}_{0},\textbf{B}_{1},...,\textbf{B}_{N-1}\}, \mathbb{C}= \{\textbf{C}_{0},\textbf{C}_{1},...,\textbf{C}_{N-1}\}, \mathbb{Q}= \{\textbf{Q}_{0},\textbf{Q}_{1},...,\textbf{Q}_{N-1}\}, \mathbb{R}= \{\textbf{R}_{0},\textbf{R}_{1},...,\textbf{R}_{N-1}\}$  and  the algorithm is given below:
\begin{algorithm}[H]
 \small
	\begin{algorithmic}[1] 
	
	\STATE Require $\mathbb{A},\mathbb{B},\mathbb{C},\mathbb{Q},\mathbb{R}$
	\STATE Initialize $\hat{\textbf{x}}_{0}$ and $\textbf{P}_{0}$  
		\FOR  {$k= 0~to~ N-1 $}
			\STATE ${\textbf{A}}_{{k}}= [\mathbb{A}]_{k+1},$ ${\textbf{B}}_{{k}}= [\mathbb{B}]_{k+1},{\textbf{C}}_{{k}}= [\mathbb{C}]_{k+1},$ ${\textbf{Q}}_{{k}}= [\mathbb{Q}]_{k+1},$ ${\textbf{R}}_{{k}}= [\mathbb{R}]_{k+1}$
		\STATE $\textbf{L}_{k}=\textbf{A}_{k}\textbf{P}_{k}\textbf{C}_{k}^{T}[\textbf{C}_{k}\textbf{P}_{k}\textbf{C}_{k}^{T}+\textbf{R}_{k}]^{-1}$ 
		\STATE Obtain $\textbf{y}_{k}$ from the sensor measurements/forward simulation of the system (17)
		\STATE $  \hat{\textbf{x}}_{k+1}=\textbf{A}_{k}\hat{\textbf{x}}_{k}+\textbf{B}_{k}\textbf{u}_{k}+\textbf{L}_{k}[\textbf{y}_{k}-\hat{\textbf{y}}_{k}]$ 
		\STATE $ \textbf{P}_{k+1}=[\textbf{A}_{k}-\textbf{L}_{k}\textbf{C}_{k}]\textbf{P}_{k}[\textbf{A}_{k}-\textbf{L}_{k}\textbf{C}_{k}]^{T}+\textbf{Q}_{k}+\textbf{L}_{k}\textbf{R}_{k}\textbf{L}_{k}^{T}$

		\ENDFOR	
	\end{algorithmic}
	\caption{(Kalman predictor)}
\end{algorithm}
Note that, the Kalman predictor algorithm consist of a forward recursion which can be combined with the forward simulation (Algorithm 2) of the system, in which we can use the estimated state for implementing the control law, i.e. $\textbf{u}_{k}=-\textbf{K}_{k}\hat{\textbf{x}}_{k}.$
\subsection{Kalman filter}
In Kalman filter, the measurements upto  $k^{th}$ time instant are used for computing the estimate $\hat{\textbf{x}}_{k},$ i.e., $l=k.$ The Kalman filter computes the estimate of the state  using the following difference equation:
 \begin{equation}
 \hat{\textbf{x}}_{k}=\textbf{A}_{k-1}\hat{\textbf{x}}_{k-1}+\textbf{B}_{k-1}\textbf{u}_{k-1}+\textbf{L}_{k}[\textbf{y}_{k}-\hat{\textbf{y}}_{k}].
 \end{equation}
 We define the estimation error vector as ${\textbf{x}}_{e_{k}}={\textbf{x}}_{k}-\hat{\textbf{x}}_{k},$ and  from (17),(36) we obtain the error dynamics as
\begin{equation}
\begin{aligned}
{\textbf{x}}_{e_{k}}&={\textbf{x}}_{k}-\hat{\textbf{x}}_{k}=\textbf{A}_{k-1}\textbf{x}_{k-1}+\textbf{B}_{k-1}\textbf{u}_{k-1}+\textbf{d}_{k-1}-\textbf{A}_{k-1}\hat{\textbf{x}}_{k-1}-\textbf{B}_{k-1}\textbf{u}_{k-1}-\textbf{L}_{k}[\textbf{C}_{k}\textbf{x}_{k}+\textbf{v}_{k}-\hat{\textbf{y}}_{k}]\\
&=\textbf{A}_{k-1}\textbf{x}_{e_{k-1}}+\textbf{d}_{k-1}-\textbf{L}_{k}\big[\textbf{C}_{k}[\textbf{A}_{k-1}\textbf{x}_{k-1}+\textbf{B}_{k-1}\textbf{u}_{k-1}+\textbf{d}_{k-1}]+\textbf{v}_{k}-\textbf{C}_{k}[\textbf{A}_{k-1}\hat{\textbf{x}}_{k-1}-\textbf{B}_{k-1}\textbf{u}_{k-1}]\big] \\
&=[\textbf{I}-\textbf{L}_{k}\textbf{C}_{k}]\big[\textbf{A}_{k-1}{\textbf{x}}_{e_{k-1}}+\textbf{d}_{k-1}\big]-\textbf{L}_{k}\textbf{v}_{k} 
 \end{aligned}
 \end{equation}
 for which we obtain the variance matrix $\textbf{P}_{k}= \textbf{V}({\textbf{x}}_{e_{k}})$ as
\begin{equation}
    \textbf{P}_{k}=[\textbf{I}-\textbf{L}_{k}\textbf{C}_{k}]\big[\textbf{A}_{k-1}\textbf{P}_{k-1}\textbf{A}_{k-1}^{T} +\textbf{Q}_{k-1}\big][\textbf{I}-\textbf{L}_{k}\textbf{C}_{k}]^{T}+\textbf{L}_{k}\textbf{R}_{k}\textbf{L}_{k}^{T}
\end{equation}
which is  the DRE for the Kalman filter. We choose the kalman gain $\textbf{L}_{k}$ to minimize the cost $J=E(\textbf{x}_{e_{k}}^{T}\textbf{x}_{e_{k}})=Trace(\textbf{P}_{k}),$ and this leads to  $\frac{\partial Trace(\textbf{P}_{k})}{\partial  \textbf{L}_{k}}=\textbf{0},$ which results in:
\begin{equation}
    \begin{aligned}
       & -2[\textbf{I}-\textbf{L}_{k}\textbf{C}_{k}]\big[\textbf{A}_{k-1}\textbf{P}_{k-1}\textbf{A}_{k-1}^{T} +\textbf{Q}_{k-1}\big]\textbf{C}_{k}^{T}+2\textbf{L}_{k}\textbf{R}_{k}=\textbf{0}\\
        &\implies \textbf{L}_{k}=[\textbf{A}_{k-1}\textbf{P}_{k-1}\textbf{A}_{k-1}^{T} +\textbf{Q}_{k-1}]\textbf{C}_{k}^{T}\big[\textbf{C}_{k}[\textbf{A}_{k-1}\textbf{P}_{k-1}\textbf{A}_{k-1}^{T} +\textbf{Q}_{k-1}]\textbf{C}_{k}^{T}+\textbf{R}_{k}\big]^{-1}.
    \end{aligned}
\end{equation}
\subsection*{Alternate derivation}
Here we can have an alternate form for Kalman filter derivation in which we separate the prediction and correction part, and this leads to a more compact form of the estimate and Riccati equation. We denote $\hat{\textbf{x}}_{k|i}$ as the estimate of $\textbf{x}_{k}$ computed using the information at $i^{th}$ time instant, i.e., we have \\
\begin{tabular}{@{}l@{\ }l}
	\hspace{.5cm}	&  $\hat{\textbf{x}}_{k|k-1}$ is the estimate of $\textbf{x}_{k}$ computed using the information at $(k-1)^{th}$ time instant. \\
	\hspace{.5cm}	 &  $\hat{\textbf{x}}_{k|k}$ is the estimate of $\textbf{x}_{k}$ computed using the information at $k^{th}$ time instant.\\
\end{tabular}\newline
Using this we can rewrite the Kalman filter equation (36) as a two-stage process:
 \begin{equation}\begin{aligned}
 &\hat{\textbf{x}}_{k|k-1}=\textbf{A}_{k-1}\hat{\textbf{x}}_{k-1|k-1}+\textbf{B}_{k-1}\textbf{u}_{k-1}\\
 &\hat{\textbf{x}}_{k|k}=\hat{\textbf{x}}_{k|k-1}+\textbf{L}_{k}[\textbf{y}_{k}-\hat{\textbf{y}}_{k|k-1}].
 \end{aligned}
 \end{equation}
and comparing this equation with (36) we obtain $\hat{\textbf{x}}_{k|k}=\hat{\textbf{x}}_{k}.$ 
Now, we can define the prediction or forecast error vector as $\textbf{x}_{e_{k|k-1}}=\textbf{x}_{k}-\hat{\textbf{x}}_{k|k-1}$ and the error dynamics:
\begin{equation}
\begin{aligned}
{\textbf{x}}_{e_{k|k-1}}&=\textbf{x}_{k}-\hat{\textbf{x}}_{k|k-1}=\textbf{A}_{k-1}\textbf{x}_{k-1}+\textbf{B}_{k-1}\textbf{u}_{k-1}+\textbf{d}_{k-1}-\textbf{A}_{k-1}\hat{\textbf{x}}_{k-1|k-1}-\textbf{B}_{k-1}\textbf{u}_{k-1}\\&=\textbf{A}_{k-1}\textbf{x}_{e_{k-1|k-1}}+\textbf{d}_{k-1}
 \end{aligned}
 \end{equation}
 for which we obtain the variance matrix $\textbf{P}_{k|k-1}= \textbf{V}({\textbf{x}}_{e_{k|k-1}})$ as
\begin{equation}
    \textbf{P}_{k|k-1}=\textbf{A}_{k-1}\textbf{P}_{k-1|k-1}\textbf{A}_{k-1}^{T} +\textbf{Q}_{k-1}.
\end{equation}
Similarly, the estimation error is defined as $\textbf{x}_{e_{k|k}}=\textbf{x}_{k}-\hat{\textbf{x}}_{k|k}$ and the error dynamics becomes
\begin{equation}
\begin{aligned}
{\textbf{x}}_{e_{k|k}}&={\textbf{x}}_{k}-\hat{\textbf{x}}_{k|k}={\textbf{x}}_{k}-\hat{\textbf{x}}_{k|k-1}-\textbf{L}_{k}[\textbf{C}_{k}\textbf{x}_{k}+\textbf{v}_{k}-\textbf{C}_{k}\hat{\textbf{x}}_{k|k-1}]\\&=[\textbf{I}-\textbf{L}_{k}\textbf{C}_{k}]\textbf{x}_{e_{k|k-1}}-\textbf{L}_{k}\textbf{v}_{k}
 \end{aligned}
 \end{equation}
 for which we obtain the variance matrix $\textbf{P}_{k|k}= \textbf{V}({\textbf{x}}_{e_{k|k}})$ as
 \begin{equation}
    \textbf{P}_{k|k}=[\textbf{I}-\textbf{L}_{k}\textbf{C}_{k}]\textbf{P}_{k|k-1}[\textbf{I}-\textbf{L}_{k}\textbf{C}_{k}]^{T}+\textbf{L}_{k}\textbf{R}_{k}\textbf{L}_{k}^{T}
\end{equation}
which is equivalent to (38) with $\textbf{P}_{k}= \textbf{P}_{k|k}.$ Then we obtain the Kalman gain from  $\frac{\partial Trace(\textbf{P}_{k|k})}{\partial  \textbf{L}_{k}}=\textbf{0},$ which gives
\begin{equation}
    \begin{aligned}
       & -2[\textbf{I}-\textbf{L}_{k}\textbf{C}_{k}]\textbf{P}_{k|k-1}\textbf{C}_{k}^{T}+2\textbf{L}_{k}\textbf{R}_{k}=\textbf{0}\\
        &\implies \textbf{L}_{k}=\textbf{P}_{k|k-1}\textbf{C}_{k}^{T}[\textbf{C}_{k}\textbf{P}_{k|k-1}\textbf{C}_{k}^{T}+\textbf{R}_{k}]^{-1}.
    \end{aligned}
\end{equation}
We can also obtain the Kalman filter equations from the conditional expectation equation by defining $\hat{\textbf{x}}_{k|k}=\textbf{E}(\textbf{x}_{k}|\textbf{y}_{k})$ and the derivation is left as an exercise.
 In Kalman filter the estimate of the state is obtained as a weighted sum of the prediction term and correction term. Also, from (44) we obatin
\begin{equation}
\begin{aligned}
\textbf{V}({\textbf{x}}_{e_{k|k}})&=\textbf{P}_{k|k-1}-2\textbf{L}_{k}\textbf{C}_{k}\textbf{P}_{k|k-1}+\textbf{L}_{k}[\textbf{C}_{k}\textbf{P}_{k|k-1}    \textbf{C}_{k}^{T}+\textbf{R}_{k}]\textbf{L}_{k}^{T} \hspace{1cm} sub. \hspace{.1cm}  \textbf{L}_{k} \hspace{.1cm}  from \hspace{.1cm} (45)\\
&=\textbf{P}_{k|k-1}-\textbf{P}_{k|k-1}\textbf{C}_{k}^{T}[\textbf{C}_{k}\textbf{P}_{k|k-1}\textbf{C}_{k}^{T}+\textbf{R}_{k}]^{-1}\textbf{C}_{k}\textbf{P}_{k|k-1}.
 \end{aligned}
 \end{equation}
 For scalar systems, we have $V(x_{e_{k|k}})=P_{k|k-1}-\frac{P_{k|k-1}^{2}C_{k}^{2}}{C_{k}^{2}P_{k|k-1}+R_{k}}$ and $V(x_{e_{k|k-1}})=P_{k|k-1}$
  which gives $V(x_{e_{k|k}})\leq V(x_{e_{k|k-1}}),$ and this is true for higher dimensional cases also. Therefore, by using the correction term (or by using the measurement information) the variance in the estimation error can be reduced, and this is the central idea of the Kalman filter. For the Kalman filter algorithm, we define the sets  $\mathbb{A}= \{\textbf{A}_{0},\textbf{A}_{1},...,\textbf{A}_{N-1} \}, \mathbb{B}= \{\textbf{B}_{0},\textbf{B}_{1},...,\textbf{B}_{N-1}\}, \mathbb{C}= \{\textbf{C}_{1},\textbf{C}_{2},...,\textbf{C}_{N}\}, \mathbb{Q}= \{\textbf{Q}_{0},\textbf{Q}_{1},...,\textbf{Q}_{N-1}\}, \mathbb{R}= \{\textbf{R}_{1},\textbf{R}_{2},...,\textbf{R}_{N}\}$  and the algorithm is given below:
\begin{algorithm}[H]
 \small
	\begin{algorithmic}[1] 
   \STATE Require $\mathbb{A}.\mathbb{B},\mathbb{C},\mathbb{Q},\mathbb{R}$
	\STATE Initialize $\hat{\textbf{x}}_{0|0}$ and $\textbf{P}_{0|0}$  
		\FOR  {$k= 1~to~ N $}
			\STATE ${\textbf{A}}_{{k-1}}= [\mathbb{A}]_{k},$ ${\textbf{B}}_{{k-1}}= [\mathbb{B}]_{k}, {\textbf{C}}_{{k}}= [\mathbb{C}]_{k},$ ${\textbf{Q}}_{{k-1}}= [\mathbb{Q}]_{k},$ ${\textbf{R}}_{{k}}= [\mathbb{R}]_{k}$
		\STATE $\hat{\textbf{x}}_{k|k-1}=\textbf{A}_{k-1}\hat{\textbf{x}}_{k-1|k-1}+\textbf{B}_{k-1}\textbf{u}_{k-1}$ 
		\STATE $ \textbf{P}_{k|k-1}=\textbf{A}_{k-1}\textbf{P}_{k-1|k-1}\textbf{A}_{k-1}^{T} +\textbf{Q}_{k-1}$ 
		 
		\STATE $\textbf{L}_{k}=\textbf{P}_{k|k-1}\textbf{C}_{k}^{T}[\textbf{C}_{k}\textbf{P}_{k|k-1}\textbf{C}_{k}^{T}+\textbf{R}_{k}]^{-1}$
		\STATE Obtain $\textbf{y}_{k}$ from the sensor measurements or forward simulation of the system (17)
		\STATE $\hat{\textbf{x}}_{k|k}=\hat{\textbf{x}}_{k|k-1}+\textbf{L}_{k}[\textbf{y}_{k}-\hat{\textbf{y}}_{k|k-1}]$
		
		\STATE $ \textbf{P}_{k|k}=[\textbf{I}-\textbf{L}_{k}\textbf{C}_{k}]\textbf{P}_{k|k-1}[\textbf{I}-\textbf{L}_{k}\textbf{C}_{k}]^{T}+\textbf{L}_{k}\textbf{R}_{k}\textbf{L}_{k}^{T}$
		
		\ENDFOR
	\end{algorithmic}
	\caption{(Kalman filter)}
\end{algorithm}
\subsection{Kalman smoother (Rauch-Tung-Striebel  smoother)}
In Kalman smoother, the measurements upto time instant $N$ are used for computing the estimate $\hat{\textbf{x}}_{k},$ i.e., $l=N.$ Kalman smoother is used for post-processing of the states, and used for applications in which the real-time estimation of the state is not required.  Kalman smoother algorithm consists of a forward recursion (filtering) followed by a backward recursion (smoothing), i.e., the algorithm consists of two stages:\\
\begin{tabular}{@{}l@{\ }l}
	\hspace{.5cm}	i. & Forward pass: in which we estimate the state $\textbf{x}_{k}$ using the measurement information upto time instant $k,$ i.e.,  \\& we  have a filtering problem  in which we compute $\hat{\textbf{x}}_{k|k}$ and ${\textbf{P}}_{k|k}$ using a forward recursion (Kalman Filter).\\
	\hspace{.5cm}	ii. &  Backward pass: in which we improve  the estimate $\hat{\textbf{x}}_{k|k}$ using the  measurements: $\{\textbf{y}_{k+1},...,\textbf{y}_{N} \},$ 
  i.e., \\& we  compute $\hat{\textbf{x}}_{k|N}$ and ${\textbf{P}}_{k|N}$ using a backward recursion.\\
\end{tabular}\newline
The first stage is already discussed in section 3.2 and here we will discuss the second stage, i.e., backward pass algorithm. We have using the conditional expectation equation (25) the expected value of $\textbf{x}_{k}$ given $\textbf{x}_{k+1}$ can be obtained:
\begin{equation}
\begin{aligned}
    \hat{\textbf{x}}_{k|k+1}& =\textbf{E}(\textbf{x}_{k}|\textbf{x}_{k+1})=  \textbf{E}(\textbf{x}_{k})+ \textbf{L} [\textbf{x}_{k+1}-\textbf{E}(\textbf{x}_{k+1})]\\
    &=\hat{\textbf{x}}_{k|k}+\textbf{L}_{s_k}[\textbf{x}_{k+1}-\hat{\textbf{x}}_{k+1|k}]
    \end{aligned}
\end{equation}
where $\textbf{L}=\textbf{L}_{s_k}\in \mathbb{R}^{n\times n}$ is the smoother gain. We define the  error  $\textbf{x}_{e_{k|k+1}}=\textbf{x}_{k}-\hat{\textbf{x}}_{k|k+1},$ and the error dynamics:
\begin{equation}
\begin{aligned}
{\textbf{x}}_{e_{k|k+1}}&=\textbf{x}_{k}-\hat{\textbf{x}}_{k|k+1}={\textbf{x}}_{k}-\hat{\textbf{x}}_{k|k}-\textbf{L}_{s_k}[\textbf{A}_{k}\textbf{x}_{k}+\textbf{B}_{k}\textbf{u}_{k}+\textbf{d}_{k}-\textbf{A}_{k}\hat{\textbf{x}}_{k|k}-\textbf{B}_{k}\textbf{u}_{k}]\\&=[\textbf{I}-\textbf{L}_{s_k}\textbf{A}_{k}]\textbf{x}_{e_{k|k}}-\textbf{L}_{s_k}\textbf{d}_{k}
 \end{aligned}
 \end{equation}
 for which  the variance matrix $\textbf{P}_{k|k+1}= \textbf{V}({\textbf{x}}_{e_{k|k+1}})$ becomes
\begin{equation}
    \textbf{P}_{k|k+1}=[\textbf{I}-\textbf{L}_{s_k}\textbf{A}_{k}]\textbf{P}_{k|k}[\textbf{I}-\textbf{L}_{s_k}\textbf{A}_{k}]^{T}+\textbf{L}_{s_k}\textbf{Q}_{k}\textbf{L}_{s_k}^{T}.
\end{equation}
 Then, we obtain the smoother gain from  $\frac{\partial Trace(\textbf{P}_{k|k+1})}{\partial  \textbf{L}_{s_k}}=\textbf{0},$ which results in:
\begin{equation}
    \begin{aligned}
       & -2[\textbf{I}-\textbf{L}_{s_k}\textbf{A}_{k}]\textbf{P}_{k|k}\textbf{A}_{k}^{T}+2\textbf{L}_{s_k}\textbf{Q}_{k}=\textbf{0}\\
        &\implies \textbf{L}_{s_k}=\textbf{P}_{k|k}\textbf{A}_{k}^{T}[\textbf{A}_{k}\textbf{P}_{k|k}\textbf{A}_{k}^{T}+\textbf{Q}_{k}]^{-1}=\textbf{P}_{k|k}\textbf{A}_{k}^{T}\textbf{P}_{k+1|k}^{-1}.
    \end{aligned}
\end{equation}
We can also derive the Kalman smoother gain from the conditional expectation equation, i.e., $\textbf{L}_{s_k}=\textbf{V}(\textbf{x}_{k},\textbf{x}_{k+1})\textbf{V}(\textbf{x}_{k+1})^{-1}=\textbf{P}_{k|k}\textbf{A}_{k}^{T}\textbf{P}_{k+1|k}^{-1},$ 
and verifying this is left as an exercise.
From (47) we can obtain the equation for the smoothed estimate $\hat{\textbf{x}}_{k|N}$ by using the law of iterated expectation (also known as the tower property), and
for the random vectors $\textbf{x},$ $\textbf{y}$ and $\textbf{z}$ we have the law of iterated expectation  as:
$\textbf{E}(\textbf{x}|\textbf{y})= \textbf{E}(\textbf{E}(\textbf{x}|\textbf{z})| \textbf{y}).$
Now, for $\textbf{x}=\textbf{x}_{k},\textbf{y}=\textbf{y}_{N},\textbf{z}=\textbf{x}_{k+1},$ we obtain 
\begin{equation}
\begin{aligned}
  \hat{\textbf{x}}_{k|N}= \textbf{E}(\textbf{x}_{k}|\textbf{y}_{N}) &=  \textbf{E}(\textbf{E}(\textbf{x}_{k}|\textbf{x}_{k+1})| \textbf{y}_{N})=\textbf{E}(\hat{\textbf{x}}_{k|k}+\textbf{L}_{s_k}[\textbf{x}_{k+1}-\hat{\textbf{x}}_{k+1|k}]  | \textbf{y}_{N})\\&=\hat{\textbf{x}}_{k|k}+\textbf{L}_{s_k}[\hat{\textbf{x}}_{k+1|N}-\hat{\textbf{x}}_{k+1|k}].
 \end{aligned}
 \end{equation}
The  error vector for the smoothed estimate is defined as $\textbf{x}_{e_{k|N}}=\textbf{x}_{k}-\hat{\textbf{x}}_{k|N},$ and the error dynamics becomes
\begin{equation}
\begin{aligned}
{\textbf{x}}_{e_{k|N}}&=\textbf{x}_{k}-\hat{\textbf{x}}_{k|N}={\textbf{x}}_{k}-\hat{\textbf{x}}_{k|k}-\textbf{L}_{s_k}[\hat{\textbf{x}}_{k+1|N}-\hat{\textbf{x}}_{k+1|k} +\textbf{x}_{k+1}-\textbf{x}_{k+1}]\\&={\textbf{x}}_{e_{k|k}}+\textbf{L}_{s_k}[\textbf{x}_{e_{k+1|N}}-\textbf{x}_{e_{k+1|k}}]
 \end{aligned}
 \end{equation}
 for which we obtain the variance matrix $\textbf{P}_{k|N}= \textbf{V}({\textbf{x}}_{e_{k|N}})$ as
\begin{equation}
    \textbf{P}_{k|N}=\textbf{P}_{k|k}+ \textbf{L}_{s_k}[\textbf{P}_{k+1|N}-\textbf{P}_{k+1|k}]\textbf{L}_{s_k}^{T}.
\end{equation}
The algorithm for the Kalman smoother is given below.
\begin{algorithm}[H]
 \small
	\begin{algorithmic}[1]

	\STATE Compute and store $
\hat{\textbf{x}}_{k|k-1},$
$\hat{\textbf{x}}_{k|k},$ $\textbf{P}_{k|k-1},$ $\textbf{P}_{k|k}$ for $k=1,2,...,N$ using Kalman filter (Algorithm 4)

		\FOR  {$k= N-1~to~ 0 $}
			\STATE ${\textbf{A}}_{{k}}= [\mathbb{A}]_{k+1}$  
		\STATE $\textbf{L}_{s_k}=\textbf{P}_{k|k}\textbf{A}_{k}^{T}\textbf{P}_{k+1|k}^{-1}$ 
		\STATE $\hat{\textbf{x}}_{k|N}=\hat{\textbf{x}}_{k|k}+\textbf{L}_{s_k}[\hat{\textbf{x}}_{k+1|N}-\hat{\textbf{x}}_{k+1|k}]$ 
		\STATE $\textbf{P}_{k|N}=\textbf{P}_{k|k}+ \textbf{L}_{s_k}[\textbf{P}_{k+1|N}-\textbf{P}_{k+1|k}]\textbf{L}_{s_k}^{T}$
		\ENDFOR
	
	\end{algorithmic}
	\caption{(Kalman smoother)}
\end{algorithm}
\subsection{Kalman estimator in steady state}
As in the LQR case, we can analyze the steady-state behavior of the Kalman estimator, in which we are interested in the convergence of the optimal estimator gain matrix and the estimator Riccati matrix.  Here we are discussing the steady-state analysis for the Kalman predictor, and similar results can be obtained for the Kalman filter and smoother.  For LTI systems,
if $[\textbf{A},\textbf{C}]$ is observable and $\textbf{Q}>0,$ the Kalman predictor DRE (29) with $\textbf{P}_{0} > 0$ converges to a unique positive definite solution $\textbf{P}$ of the Algebraic Riccati Equation (ARE):
\begin{equation}
    \textbf{P}=[\textbf{A}-\textbf{L}\textbf{C}]\textbf{P}[\textbf{A}-\textbf{L}\textbf{C}]^{T}+\textbf{Q}+\textbf{L}\textbf{R}\textbf{L}^{T}
\end{equation}
which results in the unique estimator gain: 
\begin{equation}
  \textbf{L}= \textbf{A}\textbf{P} \textbf{C}^{T}[ \textbf{C}\textbf{P}\textbf{C}^{T} +\textbf{R}]^{-1} 
\end{equation}
such that all the eigenvalues of $\textbf{A}-\textbf{L}\textbf{C}$ lies inside the unit disk. 
\par One important observation is that for Kalman estimators, the transient period for both the states and estimator gains starts from $k=0,$ i.e., during the transient period of the states, the optimal estimator gains will be time-varying. For a better understanding of this idea, we consider the estimation problem for the LTI system  with dynamics:  
\begin{equation}
    \textbf{A}=\left[\begin{matrix} 0.5 &0 \\-1 & 1.5
\end{matrix}\right] \hspace{1cm} \textbf{B}=\left[\begin{matrix} 0.5 \\0.1
\end{matrix}\right] \hspace{1cm} \textbf{C}=\left[\begin{matrix} 1 &0.5
\end{matrix}\right]
\end{equation}
and we choose simulation parameters as $\textbf{P}_{0}=\textbf{I}_{2}, \textbf{Q}=\textbf{I}_{2}, \textbf{R}=1$ and $\hat{\textbf{x}}_{0}=\left[\begin{matrix} 10 & 5
\end{matrix}\right]^{T},$ the vectors $\textbf{x}_{0},\textbf{d}_{k}, \textbf{v}_{k}$ are chosen as Gaussian random vectors, i.e. $\textbf{x}_{0}=\hat{\textbf{x}}_{0}+2.5\textbf{g}_{2}, \textbf{d}_{k}=0.25\textbf{g}_{2},\textbf{v}_{k}=0.25\textbf{g}_{1},$ where $\textbf{g}_{n}$ is the $n-$ dimensional Gaussian random vector, and the control input is chosen as $\textbf{u}_{k}=-\textbf{K}\textbf{x}_{k}$ with $\textbf{K}=\left[\begin{matrix} 2.73 &-2.75
\end{matrix}\right]$. The simulation response for the Kalman predictor,  filter and smoother are shown in Fig. 4 which contains the plots for the estimated state, elements of the gain matrix $\textbf{L}_{k}$ and variance matrix $\textbf{P}_{k}.$ From the plot of variance matrix elements we can observe that, among the three estimators, the variance is minimum for the Kalman smoother and maximum for the Kalman predictor. This indicates that the Kalman smoother gives a more reliable estimate of the states. But, the Kalman smoother cannot estimate the states in real time, i.e., there is a delay associated with the Kalman smoother estimate, hence not suitable for implementing state feedback.
In this example, the optimal estimator gain and the Riccati matrix converges to a fixed matrix, since the system is observable and there is no need of computing them once they converges. In general, we can use the fixed gain $\textbf{L}$ obtained from (55), for ${k}>{k}_\textbf{L},$ by which the online computation can be reduced.
\begin{figure} [H]
	
	\begin{center}
		\includegraphics [scale=.44] {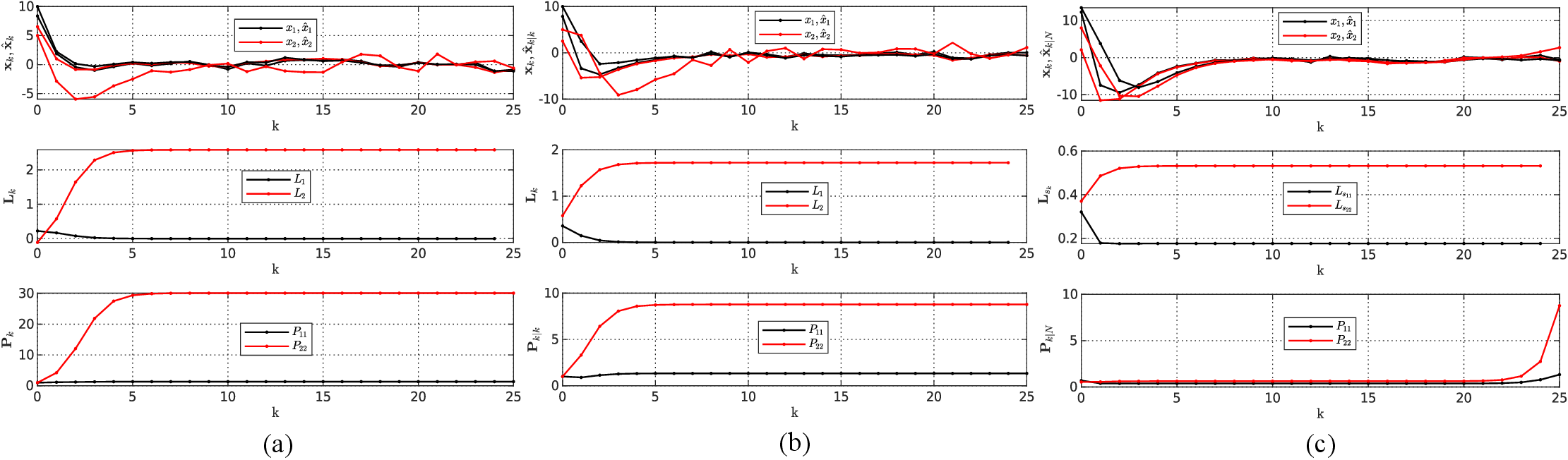}
		\caption{{\footnotesize  Simulation response of LTI system with:  (a) Kalman predictor, (b) Kalman filter, (c) Kalman smoother.}}
	\end{center}
\end{figure} 
\section{Further Reading}
This technical note attempts to discuss the basic theory  of the linear quadratic regulator and Kalman filter, which are originally proposed by R. E. Kalman in \cite{bKL,bKF} respectively. The initial contributions in the area of LQR and KF include \cite{bPA,bHS}. The well-known books that discuss the LQR problem are \cite{bDP,bBH}, and the Kalman filter and estimators are discussed in the books \cite{bBH,bAG,bGA,bFL}.
Here we have considered the discrete-time linear system model for the derivation of LQR and KF, and a continuous counterpart of these derivations can be found in \cite{bKL,bFL}.  For more advanced or related topics on LQR and KF such as: Constrained LQR (CLQR), Model Predictive Control (MPC), Extended Kalman Filter (EKF), Unscented Kalman Filter (UKF), etc., one can refer to \cite{bGA,bFL,bBM}. For more details on linear systems theory and the basic concepts such as controllability, observability, eigenvalue placement problem, state feedback controller and Luenberger observer one can refer to \cite{bCC}, and for more information on the Mayne-Murdoch formula refer to \cite{bMM}. Finally, the control approach that uses an LQR based state feedback control, in which the states are estimated using the Kalman filter is known as the Linear Quadratic Gaussian (LQG) control, i.e., an alternate, in fact, a more optimal title for this article could be: "A note on Linear Quadratic Gaussian control".

\bibliographystyle{unsrt}

\end{document}